\documentclass[12pt]{article}

\usepackage[utf8]{inputenc}

\usepackage{amsfonts}
\usepackage{amssymb}
\usepackage{amsmath}
\usepackage{amsthm}
\usepackage[T2A]{fontenc}

\usepackage[russian]{babel}

\usepackage{amsthm}

\def \le {\leqslant}
\def \ge {\geqslant}

\theoremstyle{plain}

\newtheorem{lem}{Lemma}[section]
\newtheorem{foll}{Corollary}[section]
\newtheorem{theor}{Theorem}

\newenvironment{solve}{\begin{proof}[Proof]}{\end{proof}}
\addto\captionsrussian{}
\addto\captionsrussian{}

\begin{document}
 
\centerline{\Large \bf Attainable numbers and the Lagrange spectrum}
\vspace*{3mm}
\centerline{\Large Dmitry Gayfulin\footnote{ Research is supported by RNF grant No. 14-11-00433}}
\vspace*{3mm}

\vspace*{8mm}
\begin{abstract}
For any real number $\alpha$ define the Lagrange constant $\mu(\alpha)$ by
$$
\mu^{-1}(\alpha)=\liminf_{p\in\mathbb{Z}, q\in\mathbb{N}} |q(q\alpha-p)|.
$$

The set of all values taken by $\mu(\alpha)$ as $\alpha$ varies is called the \textit{Lagrange spectrum} $\mathbb{L}$. Irrational $\alpha$ is called attainable if the inequality 
$$
\biggl|\alpha -\frac{p}{q}\biggr|\le\frac{1}{\mu(\alpha)q^2}
$$

 holds for infinitely many integers $p$ and $q$. Malyshev in survey paper \cite{Malyshev} claimed that for any $\lambda\in\mathbb{L}$ there exists an irrational $\alpha$ such that $\mu(\alpha)=\lambda$ and $\alpha$ is attainable. In the present paper we show that Malyshev's statement is incorrect and construct a counterexample. This counterexample is given by a left endpoint of a certain gap in the Lagrange spectrum. From another hand, we prove that if $\lambda$ is not a left endpoint of some gap in the Lagrange spectrum then there exists attainable $\alpha$ with $\mu(\alpha)=\lambda$. 

In addition, we give a correct proof of a theorem announced by Dietz \cite{Dietz} which describes the structure of left endpoints of gaps in the Lagrange spectrum.

\end{abstract}

\section{Introduction}
In this paper we study the properties of the Lagrange spectrum $\mathbb{L}$ which is defined as the set of all values of Lagrange constants 
$$
\mu(\alpha)=\biggl(\liminf_{p\in\mathbb{Z}, q\in\mathbb{N}} |q(q\alpha-p)|\biggr)^{-1}
$$
for real $\alpha$.

Let $A$ denote a doubly infinite sequence of positive integers $\ldots a_{-1}, a_0, a_1,\ldots$. For an arbitrary integer $i$ define
$$
\lambda_i(A)=[a_i;a_{i-1},\ldots]+[0;a_{i+1},a_{i+2},\ldots]
$$
and
$$
L(A)=\limsup\limits_{i\to\infty}\lambda_i(A), \quad M(A)=\sup\limits_{i}\lambda_i(A).
$$
It is a well known fact that the Lagrange spectrum can also be defined as the set of $L(A)$ for all sequences of positive integers $A$. The sequence $A$ is called \textit{eventually periodic} if there exists an integer $n$ such that both sequences $a_n,a_{n+1},\ldots$ and $a_{-n},a_{-n-1},\ldots$ are periodic. Here we do not suppose that these two periods coincide.

The Lagrange spectrum is a closed set \cite{Cusick1975} with minimal point $\sqrt{5}$. All the numbers of $\mathbb{L}$ which  are less than $3$ form a discrete set. They were described by Markoff. Markoff's results are discussed in details in books \cite{Cusick} and \cite{Cassels}, see also a recent paper by Bombieri \cite{Bombieri}. The complement of $\mathbb{L}$ is a countable union of \textit{maximal gaps} of the spectrum. These maximal gaps are open intervals which contain no elements of the Lagrange spectrum and their endpoints are elements of $\mathbb{L}$. In this paper we are mostly interested in properties of left endpoints of gaps in the Lagrange spectrum. In \cite{Dietz} Dietz stated that if $(a,b)$ is a maximal gap in the Lagrange spectrum then there exists an eventually periodic doubly infinite sequence $A$ such that
$$
a=M(A)=\lambda_0(A).
$$
However it was pointed by Cusick and Flahive \cite{Cusick}, p. 63 that the proof by Dietz is incorrect. Here we give a correct proof of the statement formulated by Dietz.

 Next, we consider a left endpoint 
$$
\lambda_0=[3;3,3,2,1,\overline{1,2}]+[0;2,1,\overline{1,2}]=\frac{62976-1498\sqrt{3}}{16357}\approx 3.6914708.
$$
We prove that for any $\alpha$ such that $\mu(\alpha)=\lambda_0$ the inequality
\begin{equation}
\label{attdef}
\biggl|\alpha -\frac{p}{q}\biggr|\le\frac{1}{\mu(\alpha)q^2}
\end{equation}
does not have infinitely many solutions. If an irrational $\alpha$ is such that the inequality (\ref{attdef}) has infinitely many solutions, then, following Malyshev \cite{Malyshev}, we call $\alpha$ \textit{attainable}. Here we should note that for any $\lambda<3, \lambda\in\mathbb{L}$ there exists an attainable $\alpha$ such that $\mu(\alpha)=\lambda$. In the survey paper \cite{Malyshev} Malyshev claimed   for any $\lambda\in\mathbb{L}$ there exists an irrational $\alpha$ such that $\mu(\alpha)=\lambda$ and $\alpha$ is attainable. Thus, our result gives a counterexample to a Malyshev's statement.
 
Let us discuss an equivalent definition of attainable numbers. If the inequality (\ref{attdef}) holds then the fraction $\frac{p}{q}$ should be some convergent fraction to $\alpha$ as $\mu(\alpha)\ge\sqrt{5}$. Let $\alpha$ have the following continued fraction expansion
$$
\alpha=[a_0;a_1,a_2,\ldots,a_n,\ldots].
$$
For any positive integer $i$ define 
$$
\lambda_i(\alpha)=[a_i;a_{i+1},a_{i+2},\ldots]+[a_i;a_{i-1},a_{i-2},\ldots,a_1],
$$
Denote by $\frac{p_n}{q_n}$ the $n-$th convergent fraction to $\alpha$. As
\begin{equation}
\label{minpnqn}
\biggl|\alpha -\frac{p_n}{q_n}\biggr|=\frac{1}{\lambda_{n+1}(\alpha)q_n^2}
\end{equation}
one can easily see that
\begin{equation}
\label{lalimsup}
\limsup\lambda_i(\alpha)=\mu(\alpha).
\end{equation}
From (\ref{attdef}), (\ref{minpnqn}) and (\ref{lalimsup}) one can easily deduce that $\alpha$ is attainable if and only if
$$
\lambda_i(\alpha)\ge\mu(\alpha)
$$
for infinitely many $i$.
\section{Main results}
In our first theorem we establish a counterexample to Malyshev's statement. 
\begin{theor}
\label{main}
The quadratic irrationality $\lambda_0=[3;3,3,2,1,\overline{1,2}]+[0;2,1,\overline{1,2}]$ belongs to $\mathbb{L}$, but if $\alpha$ is such that $\mu(\alpha)=\lambda_0$ then $\alpha$ is not attainable.
\end{theor}
It is a known fact that $\lambda_0$ is a left endpoint of the gap in the Markoff spectrum $\mathbb{M}$ (see \cite{Cusick}, table 2 on p.62). As $\mathbb{L}\subset\mathbb{M}$, Theorem \ref{main} implies that $\lambda_0$ is a left endpoint of a gap in the Lagrange spectrum. Our next theorem shows that each counterexample to Malyshev's statement is also a left endpoint of some gap in the Lagrange spectrum.
\begin{theor}
\label{th2}
If $\lambda\in\mathbb{L}$ is not a left endpoint of some maximal gap in the Lagrange spectrum then there exists an attainable $\alpha$ such that $\mu(\alpha)=\lambda$. 
\end{theor}
The next theorem was formulated by Dietz \cite{Dietz}.
\begin{theor}
\label{DiTh}
If $(a,b)$ is a maximal gap in $\mathbb{L}$ then there exists an eventually periodic doubly infinite sequence $A$ such that $a=\lambda_0(A)=M(A)$. Thus $a$ can be represented by a sum of two quadratic irrationalities.  
\end{theor}
As we mentioned in the previous section, Cusick and Flahive found an error in Dietz's proof. In this paper we give a new proof of this statement.

We prove Theorem \ref{main} in Sections 3-4, Theorem \ref{DiTh} in Section 5 and Theorem \ref{th2} in Section 6.
\section{The Lagrange spectrum contains $\lambda_0$}
Define the doubly infinite sequence $A_0=\ldots,a_{-1},a_0,a_1,\ldots$
$$
A_0=(\overline{2,1},1,2,3,3^*,3,2,1,\overline{1,2}),
$$
where * denotes the element $a_0$.
\begin{lem}
\label{A0lem}
$$
\sup_{i\in\mathbb{Z}}\lambda_i(A_0)=\lambda_1(A_0)=\lambda_{-1}(A_0)=\lambda_0.
$$
\end{lem}
\begin{solve}
One can easily see that $\lambda_1(A)=\lambda_{-1}(A)=\lambda_0$. Also 
$$
\lambda_0(A)=3+2[0;3,2,1,\overline{1,2}]=\frac{246+\sqrt{3}}{69}\approx 3.59032<\lambda_1(A_0).
$$

Now consider $i\ne-1, 0$ or $1$. Then
\begin{equation}
\label{13prim}
[0;a_{i+1},a_{i+2},\ldots]\le[0;\overline{1,3}], \quad [a_i;a_{i-1},a_{i-2},\ldots]\le[2;\overline{1,3}]
\end{equation}
because $a_k\le 3$ for any integer $k$. Inequalities (\ref{13prim}) give 
$$
\lambda_i(A_0)\le 2+2[0;\overline{1,3}]=\sqrt{21}-1<3.6<\lambda_1(A_0).
$$

Thus the supremum of $\lambda_i(A_0)$ is reached only at the points $i=-1$ and $i=1$.
\end{solve}

\begin{lem}
\label{inlag}
We have $\lambda_0\in\mathbb{L}$
\end{lem}
\begin{solve}
For any positive integer $n$ let $C^n$ denote the following finite sequence
$$
C^n=((2,1)_n,1,2,3,3,3,2,1,(1,2)_n),
$$
where a subscript $n$ attached to a sequence means that the sequence is repeated $n$ times. Define the irrational number $\alpha_0$ using continued fraction expansion
$$
\alpha_0=[0;C^1,C^2,\ldots,C^n,\ldots]=[0;b_1,b_2,\ldots].
$$
We show that
\begin{equation}
\label{lisula}
\limsup\lambda_i(\alpha_0)=\lambda_0.
\end{equation}
By $i_n$ denote the index of last element in the sequence $C^n$ which is equal to $3$. Then
$$
\lambda_{i_n}(\alpha_0)=\lambda_{i_n-2}(\alpha_0)=[3;3,3,2,1,(1,2)_n,\ldots]+[0;2,1,(1,2)_n,\ldots].
$$
Therefore
$$
\lim_{n\to\infty}\lambda_{i_n}(\alpha_0)=\lambda_0=\lim_{n\to\infty}\lambda_{i_n-2}(\alpha_0).
$$
From the other hand, if $i$ is not equal to $i_n$ or $i_{n-2}$ for some $n$ then $\lambda_{i}(\alpha_0)<3.6$. One can easily deduce this fact using the argument from the proof of Lemma \ref{A0lem}. Thus (\ref{lisula}) holds and the proof is completed.
\end{solve}
\section{Every $\alpha$ such that $\mu(\alpha)=\lambda_0$ is not attainable.}
We will frequently use the following classical lemma throughout the paper.
\begin{lem}\cite{Cusick}
\label{comp}
Suppose $\alpha=[a_0;a_1,\ldots,a_n,b_1,\ldots]$ and  $\beta=[a_0;a_1,\ldots,a_n,c_1,\ldots]$, where $n\ge 0$, $a_0$ is an integer, $a_1,\ldots,a_n,b_1,b_2,\ldots,c_1,c_2,\ldots$ are positive integers with $b_1\ne c_1$. Then for $n$ odd, $\alpha>\beta$ if and only if $b_1>c_1$; for $n$ even, $\alpha>\beta$ if and only if $b_1<c_1$. Also, 
$$
|\alpha-\beta|<2^{-(n-1)}.
$$
\end{lem}
Denote $\varepsilon_n=2^{-(n-1)}$.
\begin{lem}
\label{noseq}
If the inequality 
\begin{equation}
\label{lacond}
\limsup\lambda_i(\alpha)=\lambda_0
\end{equation}
holds for some $\alpha=[a_0,a_1,\ldots,a_n,\ldots]$, then no pattern from the list (\ref{ban}) occurs in the sequence $a_1,a_2,\ldots$ infinitely many times. Also, $a_i\le 3$ for almost all $i\in\mathbb{N}$.
\begin{equation}
\label{ban}
(1,3), (3,1), (2,2,3), (3,2,2), (3,2,3), (1,2,3,2,1).
\end{equation}
\end{lem}
\begin{solve}
If $a_i\ge 4$ for infinitely many $i$'s then
$$
\limsup\lambda_i(\alpha)\ge 4.
$$
We obtain a contradiction with (\ref{lacond}). 

 It follows from the definition of $\lambda_i$ that
$$
\limsup\lambda_i([0;a_1,a_2,\ldots,a_n,\ldots])=\limsup\lambda_i([0;a_n,a_{n+1}\ldots]),
$$
for any positive integer $n$. Thus, without loss of generality, one can say that $a_i\le 3$ for any $i\in\mathbb{N}$.

Denote the infinite sequence $a_1,a_2,\ldots$ by $A$. If the pattern $(3,1)$ occurs in $A$ infinitely many times, consider an arbitrary $n$ such that $a_n=3, \ a_{n+1}=1$. In this case we have
$$
\lambda_{n}(\alpha)=[3;1,a_{n+2},\ldots]+[0;a_{n-1},\ldots,a_1].
$$
As $a_i\le 3$ for any natural $i$, we use Lemma \ref{comp} and establish the lower estimate of the first term
$$
[3;1,a_{n+2},\ldots]\ge[3;1,\overline{1,3}].
$$
The second term does not exceed the finite continued fraction having length $n-1$ which is equal to $[0;3,1,3,1,\ldots]$. The least partial quotient equals $3$ if $n$ is even and equals $1$ otherwise. Denote the finite continued fraction $[0;3,1,3,1,\ldots]$ by $\alpha_n$. As $\alpha_n$ is a convergent fraction to $[0;\overline{3,1}]$, we have
$$
\alpha_n\ge [0;\overline{3,1}]-\varepsilon_n.
$$
Since $\lim\limits_{n\to\infty}\varepsilon_n=0$, we obtain the inequality
$$
\lambda_{n}(\alpha)\ge[3;1,\overline{1,3}]+[0;\overline{3,1}]-\varepsilon_n=\frac{39+4\sqrt{21}}{15}-\varepsilon_n\approx 3.82202-\varepsilon_n>\lambda_0+10^{-4}
$$
for sufficiently large $n$. We obtain a contradiction with (\ref{lacond}). 

Using the same argument one can easily deduce that the inversed pattern $(1,3)$ does not occur in $A$ infinitely many times. Without loss of generality we can say that the patterns $(3,1)$ and $(1,3)$ do not occur in $A$.

If the pattern $(3,2,2)$ occurs in $A$ infinitely many times, consider an arbitrary $n$ such that $a_n=3, a_{n+1}=2, a_{n+2}=2$. Then we have
$$
\lambda_{n}(\alpha)=[3;2,2,a_{n+3},\ldots]+[0;a_{n-1},\ldots,a_1].
$$
Since all $a_i\le 3$ and the patterns $(3,1)$ and $(1,3)$ do not occur in $A$, we obtain the following lower estimation
$$
[3;2,2,a_{n+3},\ldots]\ge[3;2,2,\overline{3,2}], \quad [0;a_{n-1},\ldots,a_1]\ge [0;\overline{3,2}]-\varepsilon_n.
$$
Thus,
$$
\lambda_{n}(\alpha)\ge[3;2,2,\overline{3,2}]+[0;\overline{3,2}]-\varepsilon_n=\frac{39+10\sqrt{15}}{21}-\varepsilon_n\approx 3.70142-\varepsilon_n>\lambda_0+10^{-4}
$$
for sufficiently large $n$. We obtain a contradiction with (\ref{lacond}). Using the same argument one can easily deduce that the inversed pattern $(2,2,3)$ does not occur in $A$ infinitely many times.

If the pattern $(3,2,3)$ occurs in $A$ infinitely many times, consider an arbitrary $n$ such that\\ $a_n=3, a_{n+1}=2, a_{n+2}=3$. Then we have
$$
\lambda_{n}(\alpha)=[3;2,3,a_{n+3},\ldots]+[0;a_{n-1},\ldots,a_1]>[3;2,2,a_{n+3},\ldots]+[0;a_{n-1},\ldots,a_1].
$$
Now we use the estimates from the previous paragraph.

Finally, consider that  the pattern $(3,2,2)$ occurs in $A$ infinitely many times. Consider an arbitrary $n$ such that $a_{n-2}=1,a_{n-1}=2,a_n=3,a_{n+1}=2,a_{n+2}=1$. In this case we have 
$$
\lambda_{n}(\alpha)=[3;2,1,a_{n+3},\ldots]+[0;2,1,a_{n-3},\ldots,a_1].
$$

As the patterns $(3,1)$ and $(1,3)$ do not occur in $A$ and, of course, all elements in $A$ do not exceed $3$, we obtain the following lower estimates
$$
[3;2,1,a_{n+3},\ldots]\ge[3;2,1,\overline{2,1}], \quad [0;2,1, a_{n-3},\ldots,a_1]\ge [0;\overline{2,1}]-\varepsilon_n.
$$
Thus,
$$
\lambda_{i_n}(\alpha)\ge[3;2,1,\overline{2,1}]+[0;\overline{2,1}]-\varepsilon_n=2+\sqrt{3}-\varepsilon_n\approx 3.73205-\varepsilon_n,
$$
which is greater than $\lambda_0+10^{-4}$ for sufficiently large $n$. We obtain a contradiction with (\ref{lacond}). Lemma is proved.
\end{solve}

\begin{lem}
\label{mainpatt}
If an irrational $\alpha=[0,a_1,\ldots,a_n,\ldots]$ is such that
$$
\limsup\lambda_i(\alpha)=\lambda_0,
$$
then there exist $N\in\mathbb{N}$ such that for any integer $n>N$ the inequality
\begin{equation}
\label{3691}
\lambda_{n}(\alpha)>3.691
\end{equation}
holds only when $a_n$ is left or right $3$ from the pattern $(1,2,3,3,3,2,1)$ in the infinite sequence\\ $A=(a_1,a_2,\ldots,a_n,\ldots)$.
\end{lem}
\begin{solve}
Lemma \ref{noseq} implies that there exists an integer $N>100$ such that in the infinite sequence $a_{N-10}, a_{N-9},\ldots,a_N,\ldots$ all elements do not exceed $3$ and all the patterns from the list (\ref{ban}) do not occur in this sequence. We will call the patterns from the list (\ref{ban}) \textit{prohibited}. Now consider an arbitrary integer $n>N$ such that the inequality \eqref{3691} holds. One can easily see that $\varepsilon_n<10^{-25}$.

Of course, $a_{n}=3$, because if $a_{n}\le 2$, then
$$
\lambda_{n}(\alpha)\le 2+2[0;\overline{1,2}]+\varepsilon_n=2\sqrt{3}+\varepsilon_n\approx 3.464101<3.691.
$$
As $a_{n}=3$, Lemma \ref{noseq} implies that only two variants are possible: either $a_{n-1}=a_{n+1}=3$, or among the numbers $a_{n-1}$ and $a_{n+1}$ one equals $2$ and another equals $3$. Consider the first case
\begin{equation}
\label{333}
\lambda_{n}(\alpha)=3+[0;3,a_{n-2},\ldots,a_1]+[0;3,a_{n+2},\ldots].
\end{equation}
Since the sequence $A$ does not contain the patterns  $(1, 3)$ and$(3,1)$, we obtain the upper estimate of (\ref{333}), using Lemma \ref{comp}
$$
[0;3,a_{n+2},\ldots]\le[0;3,\overline{3,2}],\quad [0;3,a_{n-2},\ldots,a_1]\le[0;3,\overline{3,2}]+\varepsilon_{n}.
$$
Thus
$$
\lambda_{n}(\alpha)=3+[0;3,a_{n-2},\ldots,a_1]+[0;3,a_{n+2},\ldots]\le3+2[0;3,\overline{3,2}]+\varepsilon=\frac{33-2\sqrt{15}}{7}+\varepsilon\approx 3.60772<3.691.
$$
We obtain a contradiction with (\ref{3691}). Thus, either $a_{n-1}=3, \ a_{n+1}=2$, or $a_{n-1}=2, \ a_{n+1}=3$. We consider the first case only, the second case treated in exactly the same way using the symmetry.
As the patterns $(3,2,3)$ and $(3,2,2)$ are prohibited, we have $a_{n+2}=1$. As the pattern $(1,3)$ is prohibited, $a_{n-2}$ equals $2$ or $3$. If $a_{n-2}=2$, then $a_{n-3}=1$ because the patterns $(3,2,3)$ and $(3,2,2)$ are prohibited. We have
$$
\lambda_{n}(\alpha)=[3;3,2,1,a_{n-2},\ldots,a_1]+[0;2,1,a_{n+2},\ldots].
$$
Since the patterns $(1,3)$ and $(3,1)$ are prohibited, we estimate the first term as follows
$$
[3;3,2,1,a_{n-2},\ldots,a_1]\le[3;3,2,1,\overline{2,1}]+\varepsilon_n.
$$
The second term is estimated in a similar way
\begin{equation}
\label{2332r}
[0;2,1,a_{i_k+2},\ldots]\le[0;2,1,\overline{1,2}].
\end{equation}
Eventually
\begin{equation*}
\begin{split}
\label{2332est}
\lambda_{n}(\alpha)=3+[0;3,2,1,a_{n-2},\ldots,a_1]+[0;2,1,a_{n+2},\ldots]\le\\ \le[3;3,2,1,\overline{2,1}]+[0;2,1,\overline{1,2}]+\varepsilon_n=4-\frac{2\sqrt{3}}{11}+\varepsilon_n\approx 3.68508.
\end{split}
\end{equation*}
We obtain a contradiction with (\ref{3691}). Hence $a_{n-2}=3$. As the pattern $(1,3)$ is prohibited,  $a_{n-3}$ equals $2$ or $3$. If $a_{n-3}=3$, then
\begin{equation*}
\lambda_{n}(\alpha)=[3;3,3,3,a_{n-4},\ldots,a_1]+[0;2,1,a_{n+2},\ldots].
\end{equation*}
The first term is estimated as follows
$$
[3;3,3,3,a_{n-4},\ldots,a_1]\le[3;3,3,3,3,2,1,\overline{1,2}]+\varepsilon_n.
$$
We use the estimate (\ref{2332r}) of the second term and finally obtain
\begin{equation*}
\begin{split}
\lambda_{n}(\alpha)=[3;3,3,3,a_{n-4},\ldots,a_1]+[0;2,1,a_{n+2},\ldots]\le[3;3,3,3,3,2,1,\overline{1,2}]+[0;2,1,\overline{1,2}]+\varepsilon_n=\\
=\frac{681609-16103\sqrt{3}}{177122}+\varepsilon_n\approx 3.69078.
\end{split}
\end{equation*}
We obtain a contradiction with (\ref{3691}). Hence $a_{n-3}=2$, as the patterns $(2,2,3)$ and $(3,2,3)$ are prohibited, we have $a_{n-4}=1$, and the lemma is proved.
\end{solve}
\begin{lem}
\label{nosolv}
If an irrational $\alpha=[0,a_1,\ldots,a_n,\ldots]$ is such that
$$
\limsup\lambda_i(\alpha)=\lambda_0,
$$
then there exist $N\in\mathbb{N}$ such that $\lambda_n(\alpha)<\lambda_0$ for any $n>N$.
\end{lem}
\begin{solve}
Consider the opposite. Then there exists a growing sequence $k(j)$ such that 
$$
\lambda_{k(j)}(\alpha)\ge\lambda_0\  \ \forall j\in\mathbb{N}.
$$
Lemma \ref{mainpatt} implies that there exists an integer $J$ such that $\forall j>J$ $a_{k(j)}$ is either left or right $3$ from the pattern $(1,2,3,3,3,2,1)$ in the sequence $A=a_1,a_2,\ldots$. Consider an arbitrary $j_0>J+2$. Define $n=k(j_0)$. Let $a_n$ be a right $3$ from the pattern $(1,2,3,3,3,2,1)$, then $a_{n-4}=1, a_{n-3}=2, a_{n-2}=3,\\ a_{n-1}=3, a_{n}=3, a_{n+1}=2, a_{n+2}=1$. In this case we have
\begin{equation*}
\lambda_n(\alpha)=[3;3,3,2,1,\ldots,a_1]+[0;2,1,\ldots].
\end{equation*}
We show that for the first term the following inequality holds
\begin{equation}
\label{12per}
[3;3,3,2,1,\ldots,a_1]<[3;3,3,2,1,\overline{1,2}].
\end{equation}
Consider the opposite. We note that $[3;3,3,2,1,\ldots,a_1]$ is not a convergent fraction for a quadratic irrationality $[3;3,3,2,1,\overline{1,2}]$. In the first continued fraction $a_{k(j_0-1)}=3$, while in the second the corresponding partial quotient equals $2$ or $1$. Denote the partial quotients of the first continued fraction by $[b_0;b_1,\ldots,b_{n-1}]$ and the partial quotients of the second continued fraction by\\
$[c_0;c_n,\ldots,c_{n-1},\ldots]$. Consider the minimal index $r$ such that $b_r\ne c_r$. Of course, $r>4$. If $r$ is even using Lemma \ref{comp} we obtain $b_r>c_r$. Since $c_r=2$, we have $b_r=3$. As $r$ is minimal index, $c_{r-1}=b_{r-1}=1$, and we have the prohibited pattern $(1,3)$ in $[3;3,3,2,1,\ldots,a_1]$. We obtain a contradiction with Lemma \ref{noseq}. If $r$ is odd, Lemma \ref{comp} implies that $b_r<c_r$, but $c_r=1$ and we obtain a contradiction.

A similar argument shows that
\begin{equation}
\label{12forw}
[0;2,1,\ldots,a_1]<[0;2,1,\overline{1,2}].
\end{equation}
Using (\ref{12per}) and (\ref{12forw}) we obtain
\begin{equation*}
\lambda_n(\alpha)=[3;3,3,2,1,\ldots,a_1]+[0;2,1,\ldots]<[3;3,3,2,1,\overline{1,2}]+[0;2,1,\overline{1,2}]=\lambda_0.
\end{equation*}
We have a contradiction. If $a_n$ is a left $3$ from the pattern $(1,2,3,3,3,2,1)$, a similar argument also leads to a contradiction. Lemma is proved.
\end{solve}
Now Theorem \ref{main} immediately follows from Lemmas \ref{inlag} and \ref{nosolv}.

\section{Proof of Theorem \ref{DiTh}}
A famous Freiman's theorem \cite{Freiman2} states that $[c_f,+\infty)\in\mathbb{L}$, where 
$$
c_f=4+[0;3,2,1,1,\overline{3,1,3,1,2,1}]+[0;4,3,2,2,\overline{3,1,3,1,2,1}]\approx 4.52783.
$$
 However for our purpose it is sufficient to use a weaker result of Freiman \cite{Freiman1} and Shecker \cite{Schecker} which states that $(\sqrt{21},+\infty)\in\mathbb{L}$.
\begin{foll}
If $a$ is a left endpoint of a gap of the Lagrange spectrum and $\alpha=[0;a_1,a_2,\ldots]$ is an irrational number such that $\mu(\alpha)=a$, then there exists $N$ such that $\forall n>N\ a_n\le 4$. 
\end{foll}
Consider real numbers $\alpha=[a_0;a_1,\ldots,a_n,b_1,\ldots], \beta=[a_0;a_1,\ldots,a_n,c_1,\ldots]$ such that $b_1\ne c_1$. Lemma \ref{comp} establishes the upper estimate of $|\alpha-\beta|$. However if all partial quotients of $\alpha$ and $\beta$ are bounded, a lower estimate of $|\alpha-\beta|$ can also be established. Denote $\delta_n=5^{-2(n+2)}$.

\begin{lem}
\label{lowerest}
Suppose $\alpha$ and $\beta$ satisfy the hypothesis of Lemma \ref{comp} and all partial quotients of $\alpha$ and $\beta$ do not exceed $4$. Then 
\begin{equation*}
\delta_n<|\alpha-\beta|<\varepsilon_n.
\end{equation*} 
\end{lem}
\begin{solve}
Denote by $\frac{p_n}{q_n}$  and $\frac{p'_n}{q'_n}$ the $n$-th convergent fraction to $\alpha$ and $\beta$ respectively. One can easily see that $q_n<5^n$ and $q'_n<5^n$. As $\frac{p_n}{q_n}\ne\frac{p'_n}{q'_n}$, either $\frac{p_{n+1}}{q_{n+1}}$ or $\frac{p_{n+2}}{q_{n+2}}$ lies between $\alpha$ and $\beta$. If $\frac{p_{n+1}}{q_{n+1}}$ lies between $\alpha$ and $\beta$ then $\frac{p'_{n+2}}{q'_{n+2}}$ also lies between these two numbers. Since $\frac{p_{n+1}}{q_{n+1}}\ne\frac{p'_{n+2}}{q'_{n+2}}$ we have
$$
\biggl|\frac{p_{n+1}}{q_{n+1}}-\frac{p'_{n+2}}{q'_{n+2}}\biggr|>\frac{1}{q_{n+1}q'_{n+2}}>\delta_n.
$$
The case when $\frac{p_{n+2}}{q_{n+2}}$ lies between $\alpha$ and $\beta$ is treated in exactly the same way.
\end{solve}
\begin{lem}
\label{Dirichlet}
Let $n$ be an arbitrary even positive number. Denote $N=N(n)=(2n+1)(4^{2n+1}+1)$. If $a_1, a_2,\ldots,a_N$ is an arbitrary integer sequence of length $N$ such that $1\le a_i\le 4$ for all $1\le i\le N$, then there exist two integers $n_1, n_2$ such that $a_{n_1+i}=a_{n_2+i}$ for all $0\le i\le 2n$ and $n_1\equiv n_2 \pmod 2$.
\end{lem}
\begin{solve}
There exist only $4^{2n+1}$ distinct sequences of length $2n+1$ with elements $1,2,3,4$. Consider $4^{2n+1}+1$ sequences: $(a_1, \ldots,a_{2n+1}), (a_{2n+2}, \ldots,a_{4n+2}),\ldots, (a_{(2n+1)4^{2n+1}+1},\ldots, a_{(2n+1)4^{2n+1}+2n+2})$. Dirchlet's principle implies that there exist two coinciding sequences among them. Denote these sequences by $(a_{n_1},\ldots,a_{n_1+2n})$ and $(a_{n_2},\ldots,a_{n_2+2n})$. As $n$ is even, $n_1\equiv n_2 \pmod 2$, that finishes the proof. 
\end{solve}
\begin{lem}
\label{Surgery}
Let $\gamma=[0;c_1,c_2,\ldots,c_N,\ldots]$ be an arbitrary irrational number, not a quadratic irrationality. Consider the sequence $C_N=(c_1,c_2,\ldots,c_N)$ and define two numbers $n_1$ and $n_2$ from Lemma \ref{Dirichlet}. Define two new sequences of positive integers
\begin{equation*}
\begin{split}
C_N^1=(c_1,c_2,\ldots, c_{n_1-1},c_{n_2},c_{n_2+1},\ldots,c_N),\\
C_N^2=(c_1,c_2,\ldots, c_{n_1-1},c_{n_1},\ldots,c_{n_2-1},c_{n_1},\ldots,c_{n_2-1},c_{n_2},c_{n_2+1},\ldots,c_N).
\end{split}
\end{equation*}
Let us also define two new irrational numbers:
\begin{equation*}
\begin{split}
\gamma^1=[0;c_1,c_2,\ldots, c_{n_1-1},c_{n_2},c_{n_2+1},\ldots,c_N,c_{N+1}\ldots]=[0;C_N^1,c_{N+1},\ldots],\\
\gamma^2=[0;c_1,c_2,\ldots, c_{n_1-1},c_{n_1},\ldots,c_{n_2-1},c_{n_1},\ldots,c_{n_2-1},c_{n_2},c_{n_2+1},\ldots,c_N,\ldots]=[0;C_N^2,c_{N+1},\ldots].
\end{split}
\end{equation*}
Then $\max(\gamma^1, \gamma^2)>\gamma$.
\end{lem}
\begin{solve}
Denote by $r$ minimal positive number such that $c_{n_1+r}\ne c_{n_2+r}$. As $\gamma$ is not a quadratic irrationality, $r$ exists. Suppose that $\gamma^1<\gamma$. The first different partial quotient in these two continued fractions is $c_{n_2+r}$  for $\gamma^1$ and  $c_{n_1+r}$  for $\gamma$. Now compare $\gamma^2$ and $\gamma$. The first different partial quotient in these two continued fractions is $c_{n_1+r}$  for $\gamma^2$ and  $c_{n_2+r}$  for $\gamma$. As $n_1\equiv n_2 \pmod 2$, either $\gamma^1$ or $\gamma^2$ is greater than $\gamma$.
\end{solve}
Denote $\max(\gamma^1, \gamma^2)$ by $\gamma^+$. Denote the corresponding sequence of partial quotients ($C_N^1$ or $C_N^2$) by $C_N^+$.

If all partial quotients of $\gamma$ are bounded, one can easily obtain a quantitative version of Lemma \ref{Surgery} using Lemma \ref{lowerest}.
\begin{foll}
\label{surfoll}
Let $\gamma=[0;c_1,c_2,\ldots,c_N,\ldots]$ be an arbitrary irrational number, not a quadratic irrationality. Suppose that all partial quotients of $\gamma$ do not exceed $4$. Define $\gamma^+$ as described in the previous lemma. Then
$$
\gamma^+-\gamma>\delta_{N+r}.
$$
\end{foll}
\begin{solve}
\end{solve}
\begin{lem}
\label{smallshift}
Let $\gamma=[0;c_1,c_2,\ldots,c_N,\ldots]$ and $\gamma'=[0;c'_1,c'_2,\ldots,c'_N,\ldots]$ be two irrational numbers with partial quotients not exceeding $4$. Suppose that every sequence of partial quotients of length $2n+1$ which occurs in the sequence $(c'_1,c'_2,\ldots,c'_N,\ldots)$ infinitely many times also occurs in the sequence $(c_1,c_2,\ldots,c_N,\ldots)$ infinitely many times. Then $\mu(\gamma')<\mu(\gamma)+2\varepsilon_{n}$.
\end{lem}
\begin{solve}
Let $k(j)$ be a growing sequence of positive integers such that
\begin{equation*}
\lim\limits_{j\to\infty} \lambda_{k(j)}(\gamma')=\mu(\gamma').
\end{equation*}
As the partial quotients of $\gamma'$ are limited, there exists a subsequence $k'(j)$ such that $c'_{k'(j_1)+i}=c'_{k'(j_2)+i}$ for any natural $j_1, j_2$ and $-n\le i\le n$. Denote the sequence
$$
c'_{k'(j)-n}, c'_{k'(j)-n+1},\ldots, c'_{k'(j)},\ldots, c'_{k'(j)+n}
$$
by $D$. The index is omitted because these sequences coincide for all $j$'s. The pattern $D$ has length $2n+1$ and occurs in the sequence $(c'_1,c'_2,\ldots,c'_N,\ldots)$ infinitely many times. Hence there exists a growing sequence of indices $l(j)$ such that
\begin{equation*}
c'_{k'(j)+i}=c_{l(j)+i}
\end{equation*}
for any $j\in\mathbb{N}$ and $-n\le i\le n$. Lemma \ref{lowerest} implies that
\begin{equation}
\label{smsh1}
|\lambda_{k'(j)}(\gamma') - \lambda_{l(j)}(\gamma)|<2\varepsilon_{n}.
\end{equation}
Note that 
\begin{equation}
\label{smsh2}
\lim\limits_{j\to\infty} \lambda_{k'(j)}(\gamma')=\lim\limits_{j\to\infty} \lambda_{k(j)}(\gamma')=\mu(\gamma')
\end{equation}
and
\begin{equation}
\label{smsh3}
\limsup\limits_{j\to\infty} \lambda_{l(j)}(\gamma)\le\mu(\gamma).
\end{equation}
The statement of the lemma now follows from (\ref{smsh1}), (\ref{smsh2}) and (\ref{smsh3}).
\end{solve}
\begin{lem}
\label{case1}
Let $n$ be an arbitrary positive integer. Define $N=N(n)$ from Lemma \ref{Dirichlet} and consider an arbitrary sequence $C_N$ of length $N$. Consider an arbitrary irrational $\gamma=[0;c_1,c_2,\ldots,c_N,\ldots]$ such that the pattern $C_N$ occurs in the sequence $C=(c_1,c_2,\ldots,c_m,\ldots)$ infinitely many times. Define the infinite sequence $C^+$ as follows: replace all patterns $C_N$ in the sequence $C$ by $C_N^+$. Define $\gamma'=[0;C^+]$. Then $\mu(\gamma')<\mu(\gamma)+2\varepsilon_{n}$
\end{lem}
\begin{solve}
$C^+_N$ equals $C^1_N$ or $C^2_N$. One can easily see that in both cases the conditions of Lemma \ref{smallshift} hold.
\end{solve}
Now we are ready to prove Theorem \ref{DiTh}.
\begin{solve}
Let $(a, b)$ be a maximal gap in the Lagrange spectrum. Consider an even $n$ such that $\varepsilon_n<\frac{b-a}{2}$ and $N=N(n)$ as defined in Lemma \ref{Dirichlet}. Consider an arbitrary irrational $\gamma=[0;c_1,\ldots]$ such that $\mu(\gamma)=a$. Without loss of generality we may assume that $c_i\le3\ \forall i\in\mathbb{N}$. There exists a monotonic sequence $k(j)$ such that
\begin{equation*}
\lim\limits_{j\to\infty} \lambda_{k(j)}(\gamma)=a
\end{equation*} 
and $k(j+1)-k(j)$ tends to $+\infty$ as $j\to\infty$. For each integer $j$ consider the finite sequence 
\begin{equation*}
c_{k(j)-N}, c_{k(j)-N+1},\ldots, c_{k(j)},\ldots,c_{k(j)+N}.
\end{equation*}
As all partial quotients of $\alpha$ are bounded there exists an infinite sequence of indices $j_m$ such that
$c_{k(j_{m_1})+i}=c_{k(j_{m_2})+i}$ for all $m_1, m_2\in\mathbb{N}$ and $-N\le i\le N$. Without loss of generality we can assume that $c_{k(j_1)+i}=c_{k(j_2)+i}$ for all $j_1, j_2\in\mathbb{N}$ and $-N\le i\le N$. As the sequence $(c_{k(j)+1},c_{k(j)+2},\ldots,c_{k(j)+N})$ is independent of $j$, we denote it by $C_N$. Denote $c_{k(j)}$ by $c$ as $c_{k(j)}$ is also independent of $j$. Define an infinite continued fraction 
$$
\eta_j(\gamma)=[c_{k(j)}; c_{k(j)+1}, \ldots, c_{k(j)+N}, c_{k(j)+N+1}, \ldots]=[c_{k(j)};C_N, c_{k(j)+N+1},\ldots].
$$
Lemma \ref{Dirichlet} implies that there exist two positive integers $n_1$ and $n_2$ such that
$$
c_{k(j)+n_1+i}=c_{k(j)+n_2+i}, \ \forall j\in\mathbb{N},\ 0\le i\le n-1.
$$
Denote by $r(j)$ the minimal positive number such that $c_{k(j)+n_1+r(j)}\ne c_{k(j)+n_2+r(j)}$. Of course, $r(j)\ge n$. Now consider the two cases: 
\begin{enumerate}
\item{There exists a constant M such that $r(j)=M$ for infinitely many $j$-s.}
\item{$r(j)\to \infty$ as $j\to\infty$.}
\end{enumerate}
\subsection{Case 1}
Again, without loss of generality we may assume that $\forall j\in\mathbb{N}\ r(j)=M$. Denote the segment  
$$
(c_{k(j)+N+1},\ldots, c_{k(j+1)})
$$ 
by $B_j$ and the initial segment 
$$
(c_1,c_2\ldots,c_{k(1)-1})
$$
by $B_0$. Using the introduced notation we have
\begin{equation*}
\gamma=[0;B_0,c,C_N,B_1,c,C_N,\ldots].
\end{equation*}
Define 
\begin{equation*}
\gamma'=[0;B_0,c,C^+_N,B_1,c,C^+_N,\ldots].
\end{equation*}
The index of the element $c$ before $B_i$ in the continued fraction expansion of $\gamma'$ we denote by $l(i+1)$. Then
$$ 
\lambda_{k(j)}(\gamma)=[c,C_N,B_j,C_N,\ldots]+[0;\overleftarrow{B_{j-1}},\overleftarrow{C_N},\ldots,\overleftarrow{B_0}]
$$
and
$$ 
\lambda_{l(j)}(\gamma')=[c,C^+_N,B_j,C_N,\ldots]+[0;\overleftarrow{B_{j-1}},\overleftarrow{C^+_N},\ldots,\overleftarrow{B_0}].
$$
Here $\overleftarrow{D}$ denotes an arbitrary finite sequence $D$ in reverse order. The Collorary \ref{surfoll} implies that $[c,C^+_N,B_j,C_N,\ldots]>[c,C_N,B_j,C_N,\ldots]+\delta_{N+r}$. As the length of $B_{j-1}$ tends to infinity, we have
\begin{equation*}
\lim\limits_{j\to\infty} \lambda_{l(j)}(\gamma')\ge\lim\limits_{j\to\infty} \lambda_{k(j)}(\gamma)+\delta_{N+r}\ge a+\delta_{N+r}.
\end{equation*}
Thus, $\mu(\gamma')>a$. From the other hand, Lemma \ref{case1} implies that $\mu(\gamma')<\mu(\gamma)+2\varepsilon_{n}<b$. We obtain a contradiction with the fact that $(a,b)$ is a gap in the Lagrange spectrum.
\subsection{Case 2}
Without loss of generality we may assume that $r(j)>n_2-n_1$ for any $j$. Consider $q(j)$ and $t(j)$ such that $r(j)=(n_2-n_1)q(j)+t(j)$, where $0\le t(j)<n_2-n_1$. Then
\begin{equation*}
\begin{split}
c_{k(j)+n_1}=c_{k(j)+n_2}=c_{k(j)+2n_2-n_1}=\ldots=c_{k(j)+n_1+q(j)(n_2-n_1)}\\
c_{k(j)+n_1+1}=c_{k(j)+n_2+1}=c_{k(j)+2n_2+1-n_1}=\ldots=c_{k(j)+n_1+1+q(j)(n_2-n_1)}\\
\ldots\\
c_{k(j)+n_2-1}=c_{k(j)+2n_2-n_1-1}=c_{k(j)+3n_2-2n_1-1}=\ldots=c_{k(j)+n_2-1+q(j)(n_2-n_1)}.
\end{split}
\end{equation*}
Thus, the sequence $c_{k(j)+n_1},c_{k(j)+n_1+1} \ldots,c_{k(j)+n_2-1}$ is repeated $q(j)$ times. Denote this sequence by $P$. It is independent of $j$. We have
\begin{equation*}
[c;C_N,B_j,C_N,\ldots]=[c;c_{k(j)+1},\ldots,c_{k(j)+n_1-1},\underbrace{P,P,\ldots,P}_{q(j)\ \text{times}},\ldots].
\end{equation*}
Since $r(j)$ tends to infinity, $q(j)$ also tends to infinity, hence
\begin{equation*}
\lim_{j\to\infty}[c;C_N,B_j,C_N,\ldots]=[c;c_{k(j)+1},\ldots,c_{k(j)+n_1-1},\overline{P}].
\end{equation*}
The limit is quadratic irrationality.

The fact that $[0;\overleftarrow{B_{j-1}},\overleftarrow{C_N},\ldots,\overleftarrow{B_0}]$ tends to quadratic irrationality is proved in exactly the same way. We consider the sequence $C_N=(c_{k(j)-N}, c_{k(j)-N+1},\ldots,c_{k(j)-1})$ of length $N$, find the two integers $n_1$ and $n_2$ from Lemma \ref{Dirichlet} and set $r(j)$ to be minimal natural number such that 
$$
c_{k(j)-n_1-r(j)}\ne c_{k(j)-n_2-r(j)}.
$$
The detailed proof is omitted.
\end{solve}
\section{Proof of Theorem \ref{th2}}
Let $Q$ be the set of quadratic irrationalities. It is a well known fact that $\mathbb{L}=\overline{\mu(Q)}$. Thus, if $\lambda$ is not a left endpoint of some maximal gap in the Lagrange spectrum, there are two options:
\begin{enumerate}
\item{There exists a quadratic irrationality $\gamma$ such that $\mu(\gamma)=\lambda$.}
\item{There exists a sequence of quadratic irrationalities $\gamma_n$ such that $\lim\limits_{n\to\infty}\mu(\gamma_n)=\lambda$ and the sequence $\mu(\gamma_n)$ decreases.}
\end{enumerate}

\subsection{Case 1}
Without loss of generality we may assume that $\gamma$ has purely periodic continued fraction expansion
\begin{equation*}
\gamma=[0;\overline{P}],
\end{equation*}
where $P=(c_1,c_2,\ldots,c_n)$. There exists a natural number $1\le j\le n$ such that
\begin{equation*}
\lim\limits_{m\to\infty}\lambda_{j+mn}(\gamma)=\mu(\gamma)=[c_j;c_{j+1},\ldots,c_n,\overline{P}]+[0;c_{j-1},\ldots,c_1,\overline{P}].
\end{equation*}
Note that for any finite sequence $R$ 
$$
\mu(\gamma)=\mu([0;R,\overline{P}]).
$$
Denote $[0;R,\overline{P}]$ by $\gamma'$. Denote the length of $R$ by $l$. Then
\begin{equation*}
\lambda_{j+mn+l}(\gamma')=[c_j;c_{j+1},\ldots,c_n,\overline{P}]+[0;c_{j-1},\ldots,c_1,\underbrace{P,P,\ldots,P}_{m\ \text{times}},\overleftarrow{R}].
\end{equation*}
Consider an arbitrary sequence $R$ such that 
\begin{equation*}
[0;c_{j-1},\ldots,c_1,\overleftarrow{R}]>[0;c_{j-1},\ldots,c_1,\overline{P}].
\end{equation*}
Then
\begin{equation*}
[0;c_{j-1},\ldots,c_1,\underbrace{P,P,\ldots,P}_{2m\ \text{times}},\overleftarrow{R}]>[0;c_{j-1},\ldots,c_1,\overline{P}]
\end{equation*}
for any natural $m$. Thus
\begin{equation*}
\lambda_{j+2mn+l}(\gamma')>\mu(\gamma)=\mu(\gamma')\quad \forall m\in\mathbb{N}.
\end{equation*}
As
\begin{equation*}
\lim\limits_{m\to\infty}\lambda_{j+mn+l}(\gamma')=\mu(\gamma)=\mu(\gamma'),
\end{equation*}
the theorem is proved in the first case.
\subsection{Case 2}
Denote the periods of continued fraction expansions of $\gamma_n$ by $P_n$. Without loss of generality we assume that all elements of all periods $P_n$ do not exceed $\mu(\gamma)+1$. As in the previous case we may assume that all $\gamma_n$ are purely periodic.

The following lemma easily follows from general properties of continued fractions.
\begin{lem}
\label{2N+1}
Let $\gamma_n$ be a purely periodic quadratic irrationality with period $P_n=(c^n_1,c^n_2,\ldots,c^n_{l_n})$ having length $l_n$. Consider $1\le j\le l_n$ such that
\begin{equation*}
\lim\limits_{m\to\infty}\lambda_{j+ml_n}(\gamma_n)=\lim\limits_{m\to\infty}[c^n_j;c^n_{j+1},\ldots,c^n_{l_n},\overline{P}]+[0;c^n_{j-1},\ldots,c^n_{1},\underbrace{P,P,\ldots,P}_{m\ \text{times}}]=\mu(\gamma_n).
\end{equation*}
Given any $\varepsilon>0$ there exists $N=N(\varepsilon)$ such that for any finite or infinite sequences $R$, $S$ we have
\begin{equation*}
[c^n_j;c^n_{j+1},\ldots,c^n_{l_n},\underbrace{P,P,\ldots,P}_{N\ \text{times}},R]+[0;c^n_{j-1},\ldots,c^n_1,\underbrace{P,P,\ldots,P}_{N\ \text{times}},S]>\mu(\gamma_n)-\varepsilon_n.
\end{equation*}
\end{lem}
The following lemma shows that the lengths of the periods $P_n$ tend to infinity.
\begin{lem}
\label{lenper}
Let $\gamma_n$ be a sequence of purely periodic quadratic irrationalities with period $P_n=(c^n_1,c^n_2,\ldots,c^n_{l_n})$ having length $l_n$. Let $\lambda$ be an irrational number, not a quadratic irrationality, such that $\lim\limits_{n\to\infty}\mu(\gamma_n)=\lambda$ and the sequence $\mu(\gamma_n)$ decreases. Then $l_n\to\infty$ as $n\to\infty$. Morover, there exists an infinite subsequence $P'_n$ of the sequence $P_n$ with the following properties: the first $n$ elements of the periods $P'_n,P'_{n+1},\ldots$ coincide. 
\end{lem}
\begin{solve}
Assume that $l_n$ does not tend to infinity. Then there exists an integer $M$ such that there are infinitely many periods of length $M$. As there are only finitely many sequences of elements of length $M$ with elements bounded by $\mu(\gamma)+1$, there exists a period $P'$ of length $M$ which occurs in the sequence $P_1, P_2,\ldots,P_n,\ldots$ infinitely many times. Let $i_n$ be the set of indices such that $P_{i_n}=P'$. Then
\begin{equation*}
\lambda=\lim\limits_{n\to\infty}\mu(\gamma_n)=\lim\limits_{n\to\infty}\mu([0;\overline{P_n}])=\lim\limits_{n\to\infty}\mu([0;\overline{P_{i_n}}])=\mu([0;\overline{P'}])>\lambda.
\end{equation*}
We obtain a contradiction.

Now we prove the second statement of this lemma. As all elements of $P_n$ are bounded, for any integer $m$ there exists a sequence $c'_1,c'_2,\ldots,c'_m$ and an infinite set of indices $i_n$ such that
$$
c^{i_n}_1=c'_1, c^{i_n}_2=c'_2,\ldots, c^{i_n}_m=c'_m
$$ 
for any integer $n$. Now the construction of the subsequence $P'_n$ is clear, the detailed proof is omitted.
\end{solve}

Now we are ready to prove Theorem \ref{th2}. Consider the sequence $P'_n$ from Lemma \ref{specper} and the corresponding sequence of quadratic irrationalities $\gamma'_n=[0;\overline{P'_n}]$. Define $\varepsilon_n=(\mu(\gamma_n)-\lambda)/3$ and $N(n)=N(\varepsilon)$. Consider an irrational number
\begin{equation*}
\gamma'=[0;\underbrace{P'_1,P'_1,\ldots,P'_1}_{2N(1)+1\ \text{times}},\underbrace{P'_2,P'_2,\ldots,P'_2}_{2N(2)+1\ \text{times}},\ldots,\underbrace{P'_n,P'_n,\ldots,P'_n}_{2N(n)+1\ \text{times}},\ldots].
\end{equation*}
Lemma \ref{2N+1} implies that there exist infinitely many $j-$s such that $\lambda_j(\gamma')>\lambda$. Lemma \ref{lenper} implies that  
\begin{equation*}
\mu(\gamma')=\lambda.
\end{equation*}
The theorem is proved.

\end{document}